\documentclass[11pt]{article}
\usepackage{amscd}
\usepackage[active]{srcltx}
\usepackage{bbm}
\usepackage{amsmath, multicol}
\usepackage{amssymb, epsfig}
\usepackage{epic}
\usepackage{caption}
\usepackage{graphicx}
\usepackage{float}
\newtheorem{lm}{Lemma}[section]
\newtheorem{thm}{Theorem}[section]

\newtheorem{con}{Conjecture}[section]
\newtheorem{cor}{Corollary}[section]
\newtheorem{remark}{Remark}[section]
\begin{document}
\title{\bf Diameters of Graphs with \\
Spectral Radius at most $\frac{3}{2}\sqrt{2}$}
\author{Jingfen Lan\thanks{Tsinghua University, Beijing, China, ({\tt jflan@139.com}).}
\and
Linyuan Lu
\thanks{University of South Carolina, Columbia, SC 29208,
({\tt lu@math.sc.edu}). This author was supported in part by NSF
grant DMS 1000475. }
}
\maketitle
\begin{abstract}
  The spectral radius $\rho(G)$ of a graph $G$ is the largest
  eigenvalue of its adjacency matrix.  Woo and Neumaier 
  discovered that a connected graph $G$ with $\rho(G)\leq
  \frac{3}{2}{\sqrt{2}}$ is either a dagger, an open quipu, or a
  closed quipu.  The reverse statement is not true. Many open quipus
  and closed quipus have spectral radius greater than
  $\frac{3}{2}{\sqrt{2}}$.  In this paper we proved the following results.
For any open quipu $G$ on $n$ vertices ($n\geq 6$)
  with spectral radius less than $\frac{3}{2}{\sqrt{2}}$, its diameter
  $D(G)$ satisfies $D(G)\geq (2n-4)/3$. This bound is tight.  For any
  closed quipu $G$ on $n$ vertices ($n\geq 13$) with spectral radius less than
  $\frac{3}{2}{\sqrt{2}}$, its diameter $D(G)$ satisfies $\frac{n}{3}< D(G)\leq
 \frac{2n-2}{3}$.  The upper  bound is tight while the lower bound is
asymptotically tight.

  Let $G^{min}_{n,D}$ be a graph with minimal spectral radius among
  all connected graphs on $n$ vertices with diameter $D$.
We applied the results and found $G^{min}_{n,D}$ for some range of $D$.
 For $n\geq 13$ and $D\in [\frac{n}{2}, \frac{2n-7}{3}]$, we proved that 
$G^{min}_{n,D}$
  is the graph obtained by attaching two paths of length
  $D-\lfloor\frac{n}{2}\rfloor$ and $D-\lceil\frac{n}{2}\rceil$ to a
  pair of antipodal vertices of the even cycle $C_{2(n-D)}$. Thus we
  settled a conjecture of Cioab-van~Dam-Koolen-Lee\cite{CDK}, who
  previously proved a special case $D=\frac{n+e}{2}$ for $e=1,2,3,4$.
\end{abstract}

\section{Introduction}
The {\it spectral radius} of a graph $G$, denoted by $\rho(G)$, is the largest eigenvalue of
its adjacency matrix.
Hoffman and Smith \cite{HS, Hoffman, smith}
determined all connected graphs $G$ with $\rho(G)\leq 2$.
The graphs $G$ with $\rho(G)<2$ are
 simple Dynkin Diagrams  $A_n$, $D_n$, $E_6$, $E_7$, and $E_8$,
while the graphs $G$ with $\rho(G)=2$ are
 simple extended  Dynkin
 Diagrams $\tilde A_n$, $\tilde D_n$, $\tilde E_6$, $\tilde E_7$, and $\tilde E_8$.
Cvetkovi\'c et al.~\cite{CDG} gave a nearly complete description
of all graphs $G$ with $2 < \rho(G) \leq \sqrt{2 + \sqrt{5}}$. Their description
was completed by Brouwer and Neumaier \cite{BN}.
Wang et al. \cite{wang} studied some graphs with spectral radii close to $\frac
{3}{2}{\sqrt{2}}$.
Woo and Neumaier \cite{WN} proved that any connected graph $G$ with $\sqrt{
2 + \sqrt{5}}<\rho(G)< \frac{3}{2}{\sqrt{2}}$ is one of the following graphs.
\begin{enumerate}
\item If $G$ has maximum degree at least $4$, then $G$ is a {\it
    dagger} (i.e., a tree obtained by  attaching a path to a leaf of the  star $S_5$).

\item If $G$ is a tree with maximum degree at most $3$, then $G$ is an {\it open quipu}
(see Figure \ref{opqui}).

\begin{figure}[htb]
\begin{center}
 \psfig{height=0.13\textwidth, file=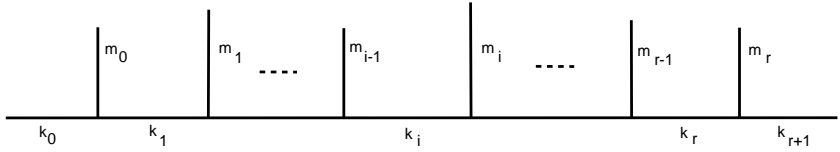}
\caption{Open Quipu $P_{(k_0,k_1,...,k_r,k_{r+1})}^{(m_0,m_1,...,m_r)}$}\label{opqui}
\end{center}
\end{figure}

\item If $G$ contains a cycle, then $G$ is a {\it closed quipu}
(see Figure \ref{cloqui}).

\begin{figure}[htb]
\begin{center}
 \psfig{height=0.2\textwidth, file=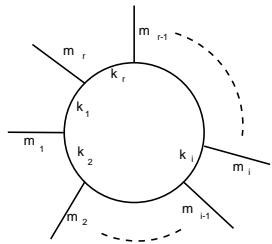}
\caption{Closed Quipu $C_{(k_1,...,k_r)}^{(m_1,...,m_r)}$}\label{cloqui}
\end{center}
\end{figure}
\end{enumerate}

No (finite) graph has spectral radius  exactly $\frac{3}{2}{\sqrt{2}}$.
The spectral radii of daggers  are always in the interval $(\sqrt{2 + \sqrt{5}}, \frac{3}{2}{\sqrt{2}})$.
 However, some open quipus (and  closed quipus) have  spectral radii greater than $\frac{3}{2}{\sqrt{2}}$.

Either an open quipu or a closed quipu can be determined by the lengths of
its internal paths and pendent paths (see Figures \ref{opqui} and
\ref{cloqui}).  Here an {\it internal path} of a graph $G$ is a path whose internal vertices
have degree $2$ and the two end vertices have degree at least $3$. An
internal path is called closed if its two end vertices coincide.
The length of an internal path is the number of its edges.
The internal path with $k$ internal vertices has length $k+1$.

 Denote by $P_{(k_0,k_1,...,k_r,k_{r+1})}^{(m_0,m_1,...,m_r)}$ the open quipu with
$r$ internal paths of lengths $k_1+1,...,k_r+1$ and $r+3$ pendent
paths of lengths $k_0, m_0,m_1,...,m_r,k_{r+1}$.  Without loss of generality,
we assume $m_0\le k_0$ and $m_r\le k_{r+1}$ through the paper. Denote by $C_{(k_1,...,k_r)}^{(m_1,...,m_r)}$ the  closed quipu with $r$ internal paths of lengths $k_1+1,...,k_r+1$ and $r$ pendent paths of lengths $m_1,...,m_r$.
Here for $1\leq i \leq r$, $k_i$ measures the number of internal vertices of the $i$-th internal path.

For convenience, a T-shape graph is viewed as
 an open quipu with $r=0$ (see Figure \ref{tt}). The graph $P_{(m+1,0,...,0,m+1)}^{(m+1,m,...,m,m+1)}$
(or  $C_{(0,...,0)}^{(m,...,m)}$) is called the {\it m-Laundry} graph
(or the {\it m-Urchin} graph) respectively.

Suppose that $G$ is a connected graph. The diameter of $G$, denoted by
$D(G)$, is the maximum distance among all pairs of vertices.
We have the following theorems.

\begin{thm} \label{thm1}
Suppose that $T$ is an open quipu on $n$ vertices $(n\geq 6)$ with
 $\rho(T)<\frac{3}{2}\sqrt{2}$. Then  the diameter of $T$ satisfies $D(T)\geq \frac{2n-4}{3}.$
The equality holds if and only if $T=P_{(1,m-2, m)}^{(1,m)}$ (for $m\geq 2$) as shown by Figure \ref{excp}.
\end{thm}

\begin{figure}[htbp]
\begin{center}
\setlength{\unitlength}{.7cm}
\begin{picture}(8,3.5)
\multiput(0,0)(1,0){3}{\circle*{0.18}}
\multiput(4,0)(1,0){3}{\circle*{0.18}}
\put(8,0){\circle*{0.18}}
\put(1,1){\circle*{0.18}}
\multiput(5,1)(0,2){2}{\circle*{0.18}}
\put(0,0){\line(1,0){2}}
\put(4,0){\line(1,0){2}}
\put(1,0){\line(0,1){1}}
\put(5,0){\line(0,1){1}}
\dashline{0.2}(2,0)(4,0)
\dashline{0.2}(6,0)(8,0)
\dashline{0.2}(5,1)(5,3)
\put(5.1,2){$m$}\put(6.7,-.5){$m$}\put(2.3,-.5){$m-2$}
\end{picture}
\hfil
\setlength{\unitlength}{.4cm}
\begin{picture}(12,7)
\put(2.2,1){\circle{4}}
\multiput(4,1)(2,0){2}{\circle*{0.36}}
\put(10,1){\circle*{0.36}}
\put(4,1){\line(1,0){2}}
\dashline{0.4}(6,1)(10,1)
\put(7.8,1.2){$m$}
\put(.6,.8){\small $2m+3$}
\put(2.2,5){\circle{4}}
\multiput(4,5)(2,0){2}{\circle*{0.36}}
\put(10,5){\circle*{0.36}}
\put(4,5){\line(1,0){2}}
\dashline{0.4}(6,5)(10,5)
\put(7.8,5.2){$m$}\put(.6,4.8){\small $2m+5$}
\end{picture}
\end{center}
\begin{multicols}{2}
\caption{$P_{(1,m-2,m)}^{(1,m)}$}\label{excp}
\caption{$C^{(m)}_{(2m+3)}$ and $C^{(m)}_{(2m+5)}$}\label{excpcq}
\end{multicols}
\end{figure}

\begin{thm}\label{thm2} Suppose that $L$ is a closed quipu on $n$ vertices $(n\geq 13)$
with  $\rho(L)< \frac{3}{2}\sqrt{2}$. Then
the diameter of $L$ satisfies $\frac{n}{3}<D(L)\le\frac{2n-2}{3}$.
Moreover, if $L$ is neither $C^{(m)}_{(2m+3)}$ nor $C^{(m)}_{(2m+5)}$
(see Figure \ref{excpcq}),  then
$D(L)\le\frac{2n-4}{3}.$
\end{thm}

\begin{remark}\label{remark1}
  The coefficient $\frac{1}{3}$ in the lower bound for $D(L)$ in Theorem \ref{thm2}
can not be improved.
  Consider the special closed quipus
  $C_{m,2m+3,r}$ with $m\ge 2$ and even $r\ge 2$ (see Figure \ref{cmkr}). Corollary
  \ref{corela} implies $\rho\left(C_{m,2m+3,r}\right)<\frac{3}{2}\sqrt
  2$ for all $m$. It has order $n=(3m+4)r$ and diameter $D=(m+2)r$. So
  $\frac{D}{n}=\frac{m+2}{3m+4}\to \frac{1}{3}$ as $m$ goes to infinity.
 \end{remark}

In 2007, van Dam and Kooij \cite{DK} asked an interesting question {\em ``which
connected graph of order $n$ with a given diameter $D$ has minimal spectral radius}?''.
A {\em  minimizer} graph, denoted by $G_{n,D}^{min}$, is a graph  which has
the minimal spectral radius among all connected graphs of order $n$
and diameter $D$. Van Dam and Kooij \cite{DK} determined $G_{n,D}^{min}$
for $D\in\{1, 2, \lfloor n/2\rfloor, n-3, n-2, n-1\}$.  The minimizer graph $G_{n,D}^{min}$
is also determined for  $D=n-4$ (Yuan-Shao-Liu \cite{YSL}),  for $D=n-5$
(Cioab\v{a}-van~Dam-Koolen-Lee \cite{CDK}),  and  for $D=n-6, n-7, n-8$ (Lan-Lu-Shi \cite{LLS}). Note $G_{n,D}^{min}$ is not unique in general.

Cioab\v{a}-van Dam-Koolen-Lee \cite{CDK} posed the following conjecture
for $D=\frac{n+e}{2}$ and proved it for $e=1,2,3,4$.

\begin{con}[Cioab\v{a}-van Dam-Koolen-Lee \cite{CDK}]
\label{con1}
For any $e\ge 1$ and sufficiently large $n$ with $n+e$  even,
$C_{(\frac{n-e-2}{2},\frac{n-e-2}{2})}^{(\lfloor\frac{e}{2}\rfloor,\lceil\frac{e}{2}\rceil)}$ is 
the unique minimizer graph $G^{min}_{n, \frac{n+e}{2}}$. 
\end{con}

We settle this conjecture by proving the statment holds for all $n\geq 3e+14$.
It is implied by the following theorem.

\begin{thm} \label{thm3}
For $n\geq 13$ and $\frac{n}{2}\le D\leq\frac{2n-7}{3}$,
$C_{(n-D-1,n-D-1)}^{(D-\lfloor\frac{n}{2}\rfloor,D-\lceil\frac{n}{2}\rceil)}$
 is the unique minimizer graph  $G^{min}_{n, D}$. 
\end{thm}

\begin{remark}\label{remark3}
It has been observed by van Dam and Kooij \cite{DK} and was finally proved by Sun \cite{Sun}
that  $C^{(m,m)}_{(2m+2, 2m+2)}$ and $P^{(m+1,m+1)}_{(m+1,2m,m+1)}$ have
the same spectral radius (see Lemma \ref{ctr}). 
Both graphs have the same $n$ and $D$.
Thus
  $C_{(n-D-1,n-D-1)}^{(D-\lfloor\frac{n}{2}\rfloor,D-\lceil\frac{n}{2}\rceil)}$
  can not be the unique minimizer graph for $n=6m+6$ and $D=3m+3$.
For $D\geq \frac{2n-2}{3}$, Sun \cite{Sun} proved that $G^{min}_{n,D}$ is
always a tree. 
\end{remark}

This paper is  organized as follows. We present some useful lemmas in section 2 and
determine the spectral radius of a family of  special quipus in section 3. 
The proofs of Theorems \ref{thm1} and \ref{thm2} are given in section 4
while the proof of Theorem \ref{thm3} is given in the last section.

\section{Basic notation and Lemmas}
\subsection{Preliminary results}
For a vertex $v$, the neighborhood of $v$ in $G$, denoted by $N(v)$, is the set
$\{u\colon uv\in E(G) \}$. Denote by  $G-v$ the remaining graph of $G$ after deleting the vertex $v$
(and all edges incident to $v$).
Similarly, $G-u-v$ is the remaining graph of $G$ after deleting two vertices $u,v$.
We need the following basic facts  (see \cite{CDS, HS, parlett, Sun}).

\begin{lm}[\cite{Sch}]
\label{lm2.1}
 Let $G$ be a graph, $v\in V(G)$, and ${\cal C}(v)$ be the set of all cycles containing $v$.
Let $e=uv$ be an edge of $G$, and ${\cal C}(e)$ be the set of all
cycles containing $e$. Then the characteristic polynomial  $\phi(G)$
satisfies
\begin{eqnarray*}
\phi(G)&=&\lambda\phi(G-v)-\sum_{w\in N(v)}\phi(G-w-v)-2\sum_{C\in
{\cal C}(v)}\phi(G-C),\\
\phi(G)&=&\phi(G-e)-\phi(G-u-v)-2\sum_{C\in {\cal C}(e)}\phi(G-C).
\end{eqnarray*}
\end{lm}

\begin{lm}[\cite{CDS}]
\label{lm2.2} Let $G_1$ and $G_2$ be two
graphs. Then the following statements hold.
\begin{enumerate}
\item If $G_2$ is a proper subgraph of $G_1$, then $\rho(G_1)>\rho(G_2)$.
\item If $\phi_{G_2}(\lambda)>\phi_{G_1}(\lambda)$ for all $\lambda\ge{\rho(G_1)}$, then $\rho(G_2)<\rho(G_1)$.
\item If $\phi_{G_1}(\rho(G_2))<0$, then $\rho(G_1)>\rho(G_2)$.
\end{enumerate}
\end{lm}

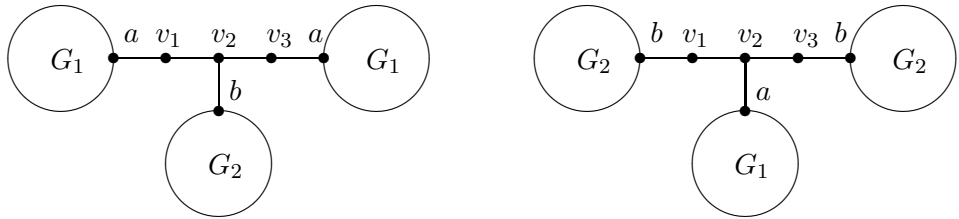
\begin{figure}[htbp]
\begin{center}
\setlength{\unitlength}{0.7cm}
\begin{picture}(20,3)
\multiput(5,0)(3,0){1}{\circle{2}}
\multiput(2,2)(6,0){2}{\circle{2}}
\multiput(5,1)(0,1){1}{\circle*{0.2}}
\multiput(3,2)(1,0){5}{\circle*{0.2}}
\multiput(5,1)(1,0){1}{\line(0,1){1}}
\multiput(3,2)(1,0){1}{\line(1,0){4}}
\put(4.8,-0.2){$G_2$}\put(5.2,1.2){$b$}\put(3.8,2.3){$v_1\quad
v_2\quad v_3$}
\multiput(3.2,2.3)(3.5,0){2}{$a$}
\multiput(1.8,1.8)(6,0){2}{$G_1$}
\multiput(15,0)(3,0){1}{\circle{2}}
\multiput(12,2)(6,0){2}{\circle{2}}
\multiput(15,1)(0,1){1}{\circle*{0.2}}
\multiput(13,2)(1,0){5}{\circle*{0.2}}
\multiput(15,1)(1,0){1}{\line(0,1){1}}
\multiput(13,2)(1,0){1}{\line(1,0){4}}
\put(14.8,-0.2){$G_1$}\put(15.2,1.2){$a$}\put(13.8,2.3){$v_1\quad
v_2\quad v_3$}
\multiput(13.2,2.3)(3.5,0){2}{$b$}
\multiput(11.8,1.8)(6,0){2}{$G_2$}
\end{picture}
\end{center}
\caption{The graphs $H_1$ and $H_2$}\label{p3}
\end{figure}

\begin{lm}[\cite{Sun}]
\label{lm2.3} Let $G_1$ and $G_2$ be two (possibly empty) graphs with $a\in V(G_1)$ and $b\in V(G_2)$, and let $H_1$ and $H_2$ be two graphs shown in Figure \ref{p3}. Then $\rho(H_1)=\rho(H_2)$.
\end{lm}

\begin{lm}[\cite{HS}]
\label{lm2.4}
Let $uv$ be an edge of a connected graph $G$ of order $n$, and denote by $G_{u,\,v}$
the graph obtained from $G$ by subdividing the edge $uv$ once, i.e.,
adding a new vertex $w$ and edges $wu,wv$ in $G-uv$. Then the
following two properties hold.
\begin{enumerate}
\item If $uv$ does not belong to an internal path of $G$ and $G\neq
C_n$, then $\rho(G_{u,\,v})>\rho(G)$.
\item If $uv$ belongs to an internal path of $G$ and $G\neq P_{(1,n-6,1)}^{(1,1)}$, then $\rho(G_{u,\,v})<\rho(G)$.
\end{enumerate}
\end{lm}

\begin{lm}[\cite{wang}]
\label{sqrt5} For any positive integer $m$, we have $$\rho(P_{(m,0,m)}^{(m,m)})<\lim_{m\to\infty}\rho(P_{(m,0,m)}^{(m,m)})=\sqrt 5.$$
\end{lm}

\begin{lm}[\cite{wang}]
\label{csh}
 For any integers $m_1,m_2\geq 1$
and  $k_1, k_2$ with $0\leq k_1 \leq k_2-2$,
we have
$$\rho(C^{(m_1,m_2)}_{(k_1,k_2)})>\rho(C^{(m_1,m_2)}_{(k_1+1,k_2-1)}).$$
\end{lm}

\begin{figure}
\begin{center}
\setlength{\unitlength}{0.7cm}
\begin{picture}(20,5)
\multiput(4,1)(7,0){3}{\circle*{0.2}}
\multiput(6,1)(7,0){3}{\circle*{0.2}}
\multiput(14,1)(1,0){4}{\circle*{0.2}}
\multiput(1,2)(2,0){2}{\circle*{0.2}}
\multiput(7,2)(2,0){2}{\circle*{0.2}}
\multiput(14,2)(3,0){2}{\circle*{0.2}}
\multiput(14,4)(3,0){2}{\circle*{0.2}}
\multiput(4,3)(2,0){2}{\circle*{0.2}}
\dashline{0.2}(4,1)(6,1)\dashline{0.2}(15,1)(16,1)
\dashline{0.2}(11,1)(13,1)\dashline{0.2}(18,1)(20,1)
\dashline{0.2}(1,2)(3,2)
\dashline{0.2}(7,2)(9,2)\dashline{0.2}(4,3)(6,3)
\dashline{0.2}(14,2)(14,4)\dashline{0.2}(17,2)(17,4)
\multiput(3,2)(3,1){2}{\line(1,-1){1}}
\multiput(3,2)(3,-1){2}{\line(1,1){1}}
\multiput(13,1)(3,0){2}{\line(1,0){2}}
\multiput(14,1)(3,0){2}{\line(0,1){1}}
\put(.9,1.9){$\underbrace{\quad\quad\quad\quad}$} \put(2,1.1){$r_1$}
\put(6.9,1.9){$\underbrace{\quad\quad\quad\quad}$} \put(8,1.1){$r_2$}
\put(3,2.5){$0$}\put(4,3.5){$1$}\put(6,3.5){$k$}\put(6.8,2.2){$k+1$}\put(3.5,.5){$2k+1$}\put(5.5,.5){$k+2$}
\put(10.9,.9){$\underbrace{\quad\quad\quad\quad}$} \put(12,.1){$r_1$}
\put(13.9,.9){$\underbrace{\quad\quad\quad\quad\quad\quad}$} \put(15.5,.1){$k$}
\put(17.9,.9){$\underbrace{\quad\quad\quad\quad}$} \put(19,.1){$r_2$}
\put(13.2,3){$r_1$}\put(17.2,3){$r_2$}
\end{picture}
\caption{$C^{(r_1-1,r_2-1)}_{(k,k)}$ and
$P_{(r_1,k-2,r_2)}^{(r_1,r_2)}$}\label{fctr}
\end{center}
\end{figure}
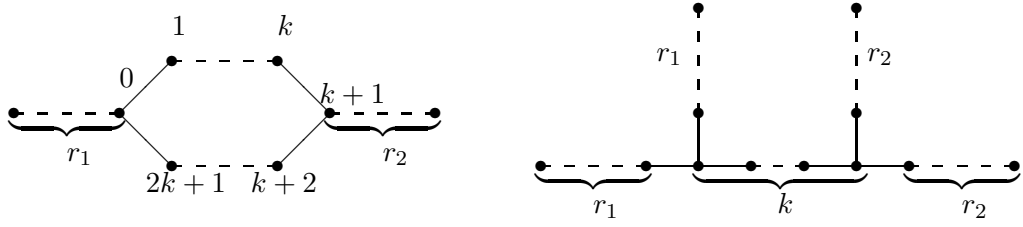

\begin{lm}[\cite{Sun}]
\label{ctr}
  For any integers $k\ge 2$, $r_1,r_2\ge 1$, we have
\begin{equation}\label{ectr}
\rho(C^{(r_1-1,r_2-1)}_{(k,k)})=\rho(P_{(r_1,k-2,r_2)}^{(r_1,r_2)}).
\end{equation}
\end{lm}
The two  graphs in Lemma \ref{ctr} are shown in Figure~\ref{fctr}.

\subsection{Our approach}
Let  $v$ be a vertex  of a graph $G$.
In \cite{LLS}, we introduced
 two functions (of $\lambda$) $p_{(G,v)}$ and $q_{(G,v)}$, which satisfy
\begin{equation}\label{eq2}
 \left(  \begin{array}[c]{c} \phi_G\\   \phi_{G-v}  \end{array}\right)
 =\left (\begin{array}{ll} 1 & 1\\ x_2 & x_1\end{array}\right )
 \left(  \begin{array}[c]{c}   p_{(G,v)}\\ q_{(G,v)} \end{array}\right).
 \end{equation}
Here $x_1$, $x_2$ are  two roots of the equation $x^2-\lambda x +1=0$.
In this paper, we always assume $\lambda\ge 2$ and $x_1\leq 1\leq x_2$.
The fact $x_1+x_2=\lambda$, $x_1x_2=1$ will be used later deliberately.
Solving $p_{(G,v)}$ and $q_{(G,v)}$, we get
\begin{equation}\label{eq3}
\left(\begin{array}[c]{c} p_{(G,v)}\\ q_{(G,v)}\end{array}\right)=\frac{1}{x_2-x_1}\left (\begin{array}{ll} -x_1 & 1\\ x_2 & -1 \end{array}\right )
\left(
  \begin{array}[c]{c}
    \phi_G\\
   \phi_{G-v}
  \end{array}
\right).
\end{equation}

For example, let $v$ be the center of the odd path $P_{2k+1}$ for $k\geq 0$.
For simplification, we 
denote $p_{(P_{2k+1},v)}$ and $q_{(P_{2k+1},v)}$ by $p_{2k+1}$ and
$q_{2k+1}$ respectively.  We have
\begin{equation}
  \label{eq:p2k+1}
 \left(\!\!
  \begin{array}[c]{c}
    p_{2k+1}\\
   q_{2k+1}
  \end{array}\!\!
\right)
=
\frac{x_2^{k+1}-x_1^{k+1}}{(x_2-x_1)^3}
\left(\!\!
\begin{array}[c]{c}
   x_2^{k-1}-2x_1^{k+1}+x_1^{k+3}\\
   x_1^{k-1}-2x_2^{k+1}+x_2^{k+3}
  \end{array}\!\!
\right).
\end{equation}

\begin{figure}[htbp]
\begin{center}
\setlength{\unitlength}{1.0cm}
\begin{picture}(2,4)
\multiput(0,.5)(1,0){2}{\circle*{0.15}}
\multiput(0,.5)(0,0){1}{\line(1,0){1}}
\multiput(1,0)(1,0){2}{\line(0,1){1}}
\multiput(1,0)(0,1){2}{\line(1,0){1}}
\put(-.1,.2){$v$}\put(1.1,.4){$v'$}\put(1.5,.3){$H$}
\put(1,-0.5){$G_0$}

\end{picture}
\hfil
\setlength{\unitlength}{1.0cm}
\begin{picture}(2,4)
\multiput(0,.5)(1,0){2}{\circle*{0.15}}
\multiput(0,1.5)(0,2){2}{\circle*{0.15}}
\multiput(0,.5)(0,0){1}{\line(1,0){1}}
\multiput(0,.5)(0,0){1}{\line(0,1){1}}
\multiput(1,0)(1,0){2}{\line(0,1){1}}
\multiput(1,0)(0,1){2}{\line(1,0){1}}
\dashline{.2}(0,1.5)(0,3.5)
\put(-.1,.2){$v$}\put(1.1,.4){$v'$}\put(1.5,.3){$H$}\put(.1,2.2){$P_m$}
\put(1,-0.5){$G_m$}
\end{picture}
\end{center}
\caption{Graph $G_m$ ($m\geq 0$).}\label{hg}
\end{figure}
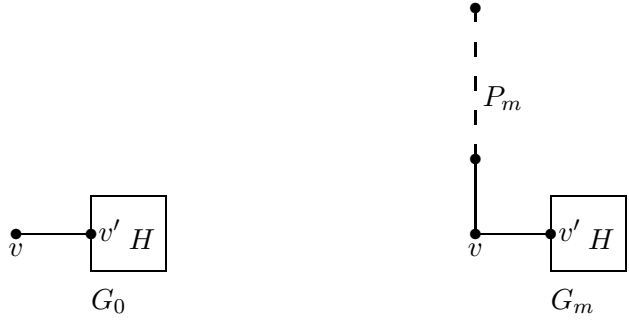

\begin{lm}\label{lmBm}
For $m\geq 0$, let $G_m$  be a graph constructed from $H$ by appending
a path $P_{m+1}$ to vertex $v'$ (see Figure \ref{hg}).
We have
$$\left(\begin{array}[c]{c} p_{(G_0,v)}\\ q_{(G_0,v)}\end{array}\right)= \left (\begin{array}{ll} x_1 & 0\\ 0 & x_2 \end{array}\right ) \left(\begin{array}[c]{c} p_{(H,v')}\\ q_{(H,v')}\end{array}\right).$$
Generally, for $m\ge 1$, we have
$$\left(\begin{array}[c]{c} p_{(G_m,v)}\\ q_{(G_m,v)}\end{array}\right)= \frac{1}{x_2-x_1}\left (\begin{array}{ll} \phi_{P_m}-x_1^{m+2} & x_1\phi_{P_{m-1}}\\ -x_2\phi_{P_{m-1}} & x_2^{m+2}-\phi_{P_m} \end{array}\right ) \left(\begin{array}[c]{c} p_{(H,v')}\\ q_{(H,v')}\end{array}\right),$$
where $\phi_{P_m}=\frac{x_2^{m+1}-x_1^{m+1}}{x_2-x_1}$.
\end{lm}

\noindent{\bf Proof }
For $m=0$, by Lemma  \ref{lm2.1}, we have
$$\left(\begin{array}[c]{c}\phi_{G_0}\\ \phi_{G_0-v} \end{array}\right)=\left (\begin{array}{ll} \lambda & -1\\ 1 & 0 \end{array}\right )\left(\begin{array}[c]{c}\phi_H\\ \phi_{H-v'} \end{array}\right).$$
Combining it with Equations \eqref{eq2} and \eqref{eq3}, we get
$$\left(\begin{array}[c]{c} p_{(G_1,v)}\\ q_{(G_1,v)}\end{array}\right)=\left (\begin{array}{ll} 1 & 1\\ x_2 & x_1 \end{array}\right )^{-1}\left (\begin{array}{ll} \lambda & -1\\ 1 & 0 \end{array}\right )\left (\begin{array}{ll} 1 & 1\\ x_2 & x_1 \end{array}\right ) \left(\begin{array}[c]{c} p_{(H,v')}\\ q_{(H,v')}\end{array}\right)\\
= \left (\begin{array}{ll} x_1 & 0\\ 0 & x_2 \end{array}\right ) \left(\begin{array}[c]{c} p_{(H,v')}\\ q_{(H,v')}\end{array}\right).$$

For $m\ge 1$, by Lemma  \ref{lm2.1}, we have
$$\left(\begin{array}[c]{c}\phi_{G_m}\\ \phi_{G_m-v} \end{array}\right)=\left (\begin{array}{ll} \phi_{P_{m+1}} & -\phi_{P_m}\\ \phi_{P_m}& 0 \end{array}\right )\left(\begin{array}[c]{c}\phi_H\\ \phi_{H-v'} \end{array}\right).$$
Similarly we get
\begin{eqnarray*}
\left(\begin{array}[c]{c} p_{(G_m,v)}\\ q_{(G_m,v)}\end{array}\right)&=&\left (\begin{array}{ll} 1 & 1\\ x_2 & x_1 \end{array}\right )^{-1}\left (\begin{array}{ll} \phi_{P_{m+1}} & -\phi_{P_m}\\ \phi_{P_m}& 0 \end{array}\right )\left (\begin{array}{ll} 1 & 1\\ x_2 & x_1 \end{array}\right ) \left(\begin{array}[c]{c} p_{(H,v')}\\ q_{(H,v')}\end{array}\right)\\
&=&\frac{1}{x_2-x_1}\left (\begin{array}{ll} \phi_{P_m}-x_1^{m+2} & x_1\phi_{P_{m-1}}\\ -x_2\phi_{P_{m-1}} & x_2^{m+2}-\phi_{P_m} \end{array}\right ) \left(\begin{array}[c]{c} p_{(H,v')}\\ q_{(H,v')}\end{array}\right).
\end{eqnarray*}
The proof is completed.
$\hfill\Box$

We define $B_m$, $d_m^{(1)}$, and $d_m^{(2)}$ as follows,
\begin{eqnarray*}
B_m&=&\frac{1}{x_2-x_1}\left (\begin{array}{ll} \phi_{P_m}-x_1^{m+2} & x_1\phi_{P_{m-1}}\\ -x_2\phi_{P_{m-1}} & x_2^{m+2}-\phi_{P_m} \end{array}\right ),\\
d_m^{(1)}&=&\phi_{P_m}-x_1^{m+2}=\frac{x_1^{m+3}-2x_1^{m+1}+x_2^{m+1}}{x_2-x_1},\\
d_m^{(2)}&=&x_2^{m+2}-\phi_{P_m}=\frac{x_2^{m+3}-2x_2^{m+1}+x_1^{m+1}}{x_2-x_1}.
\end{eqnarray*}

By a simple calculation, we have
\begin{equation}\label{m1m}
x_2^{m+2}\phi_{P_m}-d_{m+1}^{(1)}x_1^{m+1}=(x_2^{m+2}-x_1^{m+2})d_m^{(1)},
\end{equation}
\begin{equation}\label{m2m}
x_1^{m+2}\phi_{P_m}+d_{m+1}^{(2)}x_2^{m+1}=(x_2^{m+2}-x_1^{m+2})d_m^{(2)},
\end{equation}
and
\begin{equation}\label{d1d2}
 d_m^{(1)}x_2-d_m^{(2)}x_1=2\phi_{P_{m-1}}.
\end{equation}

\begin{remark}\label{remark2}
 The following equations are equivalent to each other:
\begin{eqnarray*}
&&d_m^{(2)}=\frac{2\phi_{P_{m-1}}x_1^k}{1-x_1^{k+1}},\\
&&d_m^{(2)}x_2^k-d_m^{(1)}x_1^k=2\phi_{P_{m-1}},\\
&&d_m^{(2)}=2\phi_{P_{m-1}}x_1^k+d_m^{(1)}x_1^{2k},\\
&&d_m^{(2)}=d_m^{(1)}x_1^{k-1},\\
&&d_m^{(2)}x_2^{\frac{k-1}{2}}=d_m^{(1)}x_1^{\frac{k-1}{2}}.
\end{eqnarray*}
If ``$=$'' is replaced by ``$\geq$'', then these inequalities are still equivalent to
each other.
\end{remark}
These equivalences can be proved by Equation (\ref{d1d2}).
The details are omitted.

\section{Special Quipus}

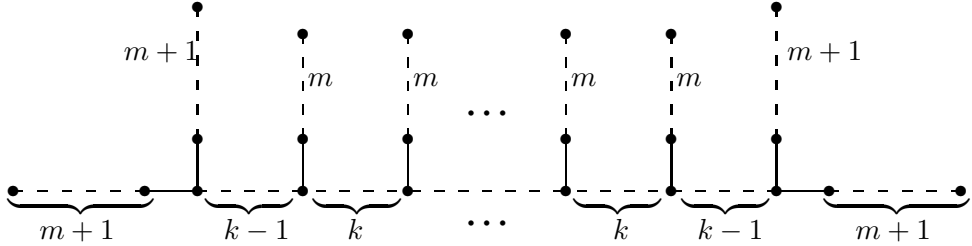
\begin{figure}[htbp]
\begin{center}
\setlength{\unitlength}{.7cm}
\begin{picture}(18,4.5)
\multiput(2,1)(2,0){3}{\circle*{0.18}}
\multiput(9,1)(2,0){3}{\circle*{0.18}}
\multiput(2,2)(2,0){3}{\circle*{0.18}}
\multiput(9,2)(2,0){3}{\circle*{0.18}}
\multiput(2,4.5)(11,0){2}{\circle*{0.18}}
\multiput(4,4)(2,0){2}{\circle*{0.18}}
\multiput(9,4)(2,0){2}{\circle*{0.18}}
\multiput(-1.5,1)(2.5,0){2}{\circle*{0.18}}
\multiput(14,1)(2.5,0){2}{\circle*{0.18}}
\multiput(2,1)(2,0){3}{\line(0,1){1}}
\multiput(9,1)(2,0){3}{\line(0,1){1}}
\multiput(1,1)(12,0){2}{\line(1,0){1}}
\dashline{0.2}(2,1)(13,1)
\dashline{0.2}(-1.5,1)(1,1)
\dashline{0.2}(14,1)(16.5,1)
\dashline{0.2}(2,2)(2,4.5)
\dashline{0.2}(13,2)(13,4.5)
\dashline{0.2}(4,2)(4,4)
\dashline{0.2}(6,2)(6,4)
\dashline{0.2}(9,2)(9,4)
\dashline{0.2}(11,2)(11,4)
\multiput(7.2,2.5)(0.3,0){3}{\circle*{0.1}}
\multiput(7.2,0.4)(0.3,0){3}{\circle*{0.1}}
\put(2.5,0.1){$k-1$}\put(4.85,0.1){$k$}\put(9.9,0.1){$k$}\put(11.5,0.1){$k-1$}
\put(-1,0.1){$m+1$}\put(14.5,0.1){$m+1$}\put(.6,3.5){$m+1$}\put(13.2,3.5){$m+1$}\put(4.1,3){$m$}\put(6.1,3){$m$}\put(9.1,3){$m$}\put(11.1,3){$m$}
\put(2.15,0.9){$\underbrace{\quad\quad\quad}$}
\put(4.2,0.9){$\underbrace{\quad\quad\quad}$}
\put(9.15,0.9){$\underbrace{\quad\quad\quad}$}
\put(11.2,0.9){$\underbrace{\quad\quad\quad}$}
\put(-1.6,0.9){$\underbrace{\quad\quad\quad\quad\quad}$}
\put(13.9,0.9){$\underbrace{\quad\quad\quad\quad\quad}$}
\end{picture}
\end{center}
\caption{A family of special trees: $P_{m,k,r}$ for $r\ge 2$ and $k\ge 1$.}\label{mkr}
\end{figure}

\begin{figure}[htbp]
\begin{center}
\setlength{\unitlength}{.7cm}
\begin{picture}(12.3,4)
\multiput(0,0)(12,0){2}{\circle*{0.18}}
\multiput(2.5,0)(1,0){3}{\circle*{0.18}}
\multiput(7.5,0)(1,0){3}{\circle*{0.18}}
\multiput(3.5,1)(5,0){2}{\circle*{0.18}}
\multiput(3.5,3)(5,0){2}{\circle*{0.18}}
\multiput(2.5,0)(5,0){2}{\line(1,0){2}}
\multiput(3.5,0)(5,0){2}{\line(0,1){1}}
\dashline{0.2}(0,0)(2.5,0)
\dashline{0.2}(4.5,0)(7.5,0)
\dashline{0.2}(9.5,0)(12,0)
\dashline{0.2}(3.5,1)(3.5,3)
\dashline{0.2}(8.5,1)(8.5,3)
\put(.5,-0.8){$m+1$}\put(10,-0.8){$m+1$}\put(2,2){$m+1$}\put(8.7,2){$m+1$}\put(5.5,-0.8){$k-2$}
\put(-.1,-.1){$\underbrace{\quad\quad\quad\quad\quad}$}
\put(9.4,-.1){$\underbrace{\quad\quad\quad\quad\quad}$}
\put(4.35,-.1){$\underbrace{\quad\quad\quad\quad\quad\quad}$}
\end{picture}
\end{center}
\caption{A family of special trees: $P_{m,k,1}$  for $k\ge 2$.}\label{mk1}
\end{figure}
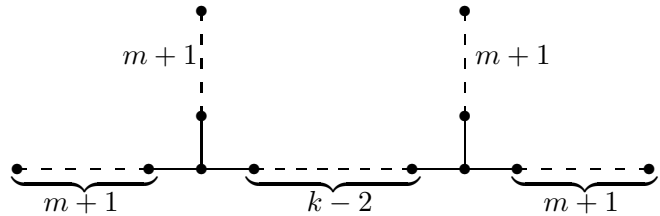

It  has been already  known that $\rho(P_{(1,n-6,1)}^{(1,1)})=2$ and $\rho(C_n)=2$ for all $n\ge 6$. This is actually a trivial case ($m=0$) of the special quipus we will show in this section.

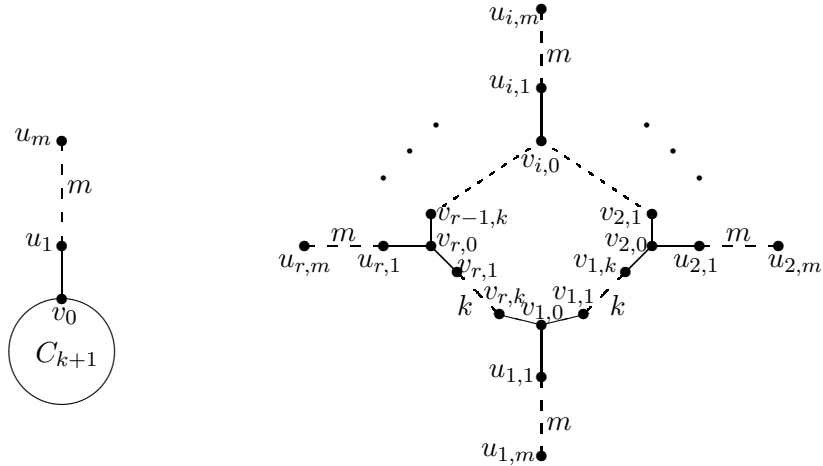
\begin{figure}[h]
\begin{center}
\setlength{\unitlength}{.7cm}
\begin{picture}(4,10)
\put(2,2){\circle{4}}
\multiput(2,3)(0,1){2}{\circle*{0.18}}
\put(2,6){\circle*{0.18}}
\put(2,3){\line(0,1){1}}
\dashline{0.2}(2,4)(2,6)
\put(2.1,5){$m$}\put(1.5,1.8){$C_{k+1}$}
\put(1.1,6){$u_m$}\put(1.3,4){$u_1$}\put(1.8,2.6){$v_0$}
\end{picture}
\hfil
\setlength{\unitlength}{.7cm}
\begin{picture}(10,8.5)
\multiput(4,0)(0,1.5){2}{\circle*{0.18}}
\multiput(4,2.5)(0,3.5){2}{\circle*{0.18}}
\multiput(3.2,2.7)(1.6,0){2}{\circle*{0.18}}
\multiput(2.4,3.5)(3.2,0){2}{\circle*{0.18}}
\multiput(1.9,4)(4.2,0){2}{\circle*{0.18}}
\multiput(1.9,4.6)(4.2,0){2}{\circle*{0.18}}
\multiput(4,7)(0,1.5){2}{\circle*{0.18}}
\multiput(-.5,4)(1.5,0){2}{\circle*{0.18}}
\multiput(7,4)(1.5,0){2}{\circle*{0.18}}
\multiput(4,1.5)(0,4.5){2}{\line(0,1){1}}
\multiput(1,4)(5.2,0){2}{\line(1,0){.9}}
\put(4,2.5){\line(4,1){.9}}
\put(4,2.5){\line(-4,1){.9}}
\multiput(1.9,4)(4.2,0){2}{\line(0,1){.6}}
\put(2.4,3.5){\line(-1,1){.55}}
\put(5.6,3.5){\line(1,1){.55}}
\dashline{0.2}(3.2,2.7)(2.4,3.5)
\dashline{0.2}(4.8,2.7)(5.6,3.5)
\dashline{0.2}(4,0)(4,1.5)
\dashline{0.2}(4,7)(4,8.5)
\dashline{0.2}(-.5,4)(1,4)
\dashline{0.2}(7,4)(8.5,4)
\dashline{0.2}(-.5,4)(1,4)
\dashline{0.1}(1.9,4.6)(4,6)
\dashline{0.1}(6.1,4.6)(4,6)
\multiput(1,5.3)(.5,.5){3}{\circle*{0.1}}
\multiput(7,5.3)(-.5,.5){3}{\circle*{0.1}}
\put(4.1,.5){$m$}\put(4.1,7.5){$m$}\put(0,4.1){$m$}\put(7.5,4.1){$m$}
\put(2.4,2.7){$k$}\put(5.3,2.7){$k$}
\put(3.6,2.7){$v_{1,0}$}\put(4.2,3){$v_{1,1}$}\put(4.6,3.6){$v_{1,k}$}
\put(2,4){$v_{r,0}$}\put(2.4,3.5){$v_{r,1}$}\put(2.9,3){$v_{r,k}$}
\put(5.2,4){$v_{2,0}$}\put(5.1,4.5){$v_{2,1}$}\put(3.6,5.5){$v_{i,0}$}\put(2,4.5){$v_{r-1,k}$}
\put(3,1.5){$u_{1,1}$}\put(2.85,0){$u_{1,m}$}
\put(6.5,3.6){$u_{2,1}$}\put(8.3,3.6){$u_{2,m}$}
\put(3,7){$u_{i,1}$}\put(3,8.3){$u_{i,m}$}
\put(.5,3.6){$u_{r,1}$}\put(-1,3.6){$u_{r,m}$}
\end{picture}
\end{center}
\caption{$C_{m,k,1}$ and $C_{m,k,r}$ ($r\ge 2$)}\label{cmkr}
\end{figure}

For integers $m,k,r\ge 1$, we respectively denote by $P_{m,k,r}$ and
$C_{m,k,r}$ the open quipu
$P_{(m+1,k-1,k,...,k,k-1,m+1)}^{m+1,m,...,m,m+1}$
($P_{(m+1,k-2,m+1)}^{(m+1,m+1)}$ for $r=1$) and the closed quipu
$C_{(k,...,k)}^{(m,...,m)}$, both of which have $r$ internal
paths. See Figure \ref{mkr}, \ref{mk1}, and \ref{cmkr}.

The quipus $P_{m,k,r}$ and $C_{m,k,r}$ play an important role in our
paper. We will show that they have the same spectral radius, which
does not depend on $r$.

For any positive $s$, we define $A^s=\left (\begin{array}{ll} x_1^s & 0\\ 0 & x_2^s \end{array}\right )$. For $k=2s+1$,
we have the following equation to use later
\begin{equation}\label{aba}
A^sB_mA^{s+1}=\frac{1}{x_2-x_1}\left (\begin{array}{ll} d_m^{(1)}x_1^k & \phi_{P_{m-1}}\\ -\phi_{P_{m-1}} & d_m^{(2)}x_2^k \end{array}\right ).
\end{equation}

\begin{lm}\label{sptree} For any integers $r, m, k\geq 1$ (except
for $r=k=1$), the spectral radius of the open quipu $P_{m,k,r}$ is the largest root $\rho_{m,k}$ of the equation $d_m^{(2)}=\frac{2\phi_{P_{m-1}}x_1^k}{1-x_1^{k+1}}$.
\end{lm}

\noindent{\bf Proof } Let $v$ be the leftmost vertex of $P_{m,k,r}$ and $s=(k-1)/2$.
 For $r\ge 2$, by Lemma \ref{lm2.1} and Lemma \ref{lmBm}, we have
\begin{eqnarray*}
\phi_{P_{m,k,r}}&=&(1, 1) \left(\begin{array}[c]{c} p_{(P_{m,k,r},v)}\\ q_{(P_{m,k,r},v)}\end{array}\right)\\
&=& (1, 1)A^{m+1}B_{m+1}A^{k-1}B_mA^k...B_mA^{k-1}B_{m+1}A^m\left(\!\!
  \begin{array}[c]{c}
    p_1\\
   q_1
  \end{array}\!\!
\right)\\
&=& \frac{1}{(x_2-x_1)^3} (d_{m+1}^{(1)}x_1^{m+1}-x_2^{m+2}\phi_{P_m}, \quad x_1^{m+2}\phi_{P_m}+d_{m+1}^{(2)}x_2^{m+1})\\
& &A^s(A^sB_mA^{s+1})^{r-1}A^s
\left( \begin{array}{ll} x_2^{m+2}\phi_{P_m}-d_{m+1}^{(1)}x_1^{m+1}\\ x_1^{m+2}\phi_{P_m}+d_{m+1}^{(2)}x_2^{m+1} \end{array}\right)\\
&=& \frac{(x_2^{m+2}-x_1^{m+2})^2}{(x_2-x_1)^3} (-d_m^{(1)}x_1^s, \quad d_m^{(2)}x_2^s)(A^sB_mA^{s+1})^{r-1}
\left( \begin{array}{ll} d_m^{(1)}x_1^s\\ d_m^{(2)}x_2^s \end{array}\right).
\end{eqnarray*}
In the last step, we applied Equations (\ref{m1m}) and (\ref{m2m}).

Now we prove that $\rho_{m,k}$ is a root of $\phi_G$.
At $\lambda=\rho_{m,k}$, by Remark \ref{remark2} we have
$$d_m^{(2)}x_2^s=d_m^{(1)}x_1^s \qquad \mbox{and} \qquad d_m^{(2)}x_2^k-\phi_{P_{m-1}}=d_m^{(1)}x_1^k+\phi_{P_{m-1}}.$$
Thus, by Equation (\ref{aba}) we get
$$(A^sB_mA^{s+1})\left( \begin{array}{ll} 1\\ 1\end{array}\right)
=\frac{1}{x_2-x_1}\left (\begin{array}{ll} d_m^{(1)}x_1^k & \phi_{P_{m-1}}\\ -\phi_{P_{m-1}} & d_m^{(2)}x_2^k \end{array}\right )\left( \begin{array}{ll} 1\\ 1\end{array}\right)
=\frac{d_m^{(1)}x_1^k+\phi_{P_{m-1}}}{x_2-x_1} \left( \begin{array}{ll} 1\\ 1\end{array}\right).$$

At the point $\lambda=\rho_{m,k}$, we have
$$
 \phi_{P_{m,k,r}}(\rho_{m,k})
= \frac{(x_2^{m+2}-x_1^{m+2})^2}{(x_2-x_1)^{r+2}}(d_m^{(1)}x_1^k+\phi_{P_{m-1}})^{r-1} (d_m^{(1)})^2x_1^{k-1} (-1,1)
\left( \begin{array}{ll} 1\\ 1\end{array}\right)= 0.
$$

It remains to prove $\phi_G(\lambda)>0$ for all $\lambda>\rho_{m,k}$.

By Remark \ref{remark2}, for $\lambda>\rho_{m,k}$, we have 
$d_m^{(2)}x_2^k-\phi_{P_{m-1}}>d_m^{(1)}x_1^k+\phi_{P_{m-1}}
$ (and $d_m^{(2)}x_2^s>d_m^{(1)}x_1^s$).
Observe that
$A^sB_mA^{s+1}$ maps the region  $\{(z_1,z_2)\colon z_2\geq z_1>0\}$
into  itself. By induction on $r$,
$(A^sB_mA^{s+1})^{r-1}$ maps the region  $\{(z_1,z_2)\colon z_2\geq z_1>0\}$
into  itself.

Let $$\left( \begin{array}{ll} z_1\\ z_2\end{array}\right)
=(A^sB_mA^{s+1})^{r-1}\left( \begin{array}{ll} d_m^{(1)}x_1^s\\ d_m^{(2)}x_2^s \end{array}\right).$$

Since $d_m^{(2)}x_2^s>d_m^{(1)}x_1^s>0$ for all $\lambda>\rho_{m,k}$, we have $z_2>z_1>0$.
Thus,
$$\phi_{P_{m,k,r}}(\rho_{m,k})=\frac{(x_2^{m+2}-x_1^{m+2})^2}{(x_2-x_1)^3} (d_m^{(2)}x_2^sz_2-d_m^{(1)}x_1^sz_1)>0.$$

For $r=1$, by the similar calculation, we have

\begin{eqnarray*}
\phi_{P_{m,k,1}}
&=& (1, 1)A^{m+1}B_{m+1}A^{k-2}B_{m+1}A^m\left(\!\!
  \begin{array}[c]{c}
    p_1\\
   q_1
  \end{array}\!\!
\right)\\
&=&\frac{(x_2^{m+2}-x_1^{m+2})^2}{(x_2-x_1)^3} (-d_m^{(1)}, d_m^{(2)})A^{k-1} \left( \begin{array}{ll} d_m^{(1)}\\ d_m^{(2)} \end{array}\right)\\
&=& \frac{(x_2^{m+2}-x_1^{m+2})^2}{(x_2-x_1)^3}\left((d_m^{(2)})^2x_2^{k-1}-(d_m^{(1)})^2x_1^{k-1}\right ).
\end{eqnarray*}

So, $\lambda=\rho(P_{m,k,1})$ is the largest root of  $d_m^{(2)}x_2^{\frac{k-1}{2}}-d_m^{(1)}x_1^{\frac{k-1}{2}}=0,$ which is equivalent to $d_m^{(2)}=\frac{2\phi_{P_{m-1}}x_1^k}{1-x_1^{k+1}}$ by Remark \ref{remark2}.

The proof of the lemma is finished.
$\hfill\Box$

\begin{lm}\label{spcqu} For $m\geq 1$, $k\geq 2$, and $r\geq 1$,
the spectral radius of the closed quipu
  $C_{m,k,r}$  is  also $\rho_{m,k}$.
\end{lm}

\noindent{\bf Proof } We observe that $C_{m,k,r}$ is a graph covering of
$C_{m,k,1}$.  The spectrum of $C_{m,k,1}$ is a subset of $C_{m,k,r}$.
The Perron-Frobenius vector of $C_{m,k,1}$ can be lifted
as the Perron-Frobenius vector of $C_{m,k,r}$. Hence,
$\rho(C_{m,k,1})= \rho(C_{m,k,r})$ for all $r\geq 2$.  By Lemmas
\ref{ctr} and \ref{sptree}, we have $\rho(C_{m,k,2})=\rho(P_{m,k,1})=\rho_{m,k}$ for
$k\ge 2$.  Hence, $\rho(C_{m,k,r})=\rho_{m,k}$ for all $r\geq 1$ and $k\geq 2$.
\hfill $\square$

\section{Quipus with spectral radii bounded by $\frac{3\sqrt 2}{2}$}

In this section, we will describe those open quipus and closed quipus
with spectral radii less than $\frac{3}{2}\sqrt 2$.

\subsection{A Lemma}
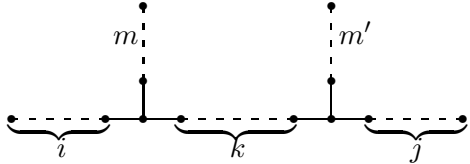
\begin{figure}[h]
\begin{center}
\setlength{\unitlength}{.5cm}
\begin{picture}(12.3,4)(0,-1)
\multiput(0,0)(12,0){2}{\circle*{0.18}}
\multiput(2.5,0)(1,0){3}{\circle*{0.18}}
\multiput(7.5,0)(1,0){3}{\circle*{0.18}}
\multiput(3.5,1)(5,0){2}{\circle*{0.18}}
\multiput(3.5,3)(5,0){2}{\circle*{0.18}}
\multiput(2.5,0)(5,0){2}{\line(1,0){2}}
\multiput(3.5,0)(5,0){2}{\line(0,1){1}}
\dashline{0.2}(0,0)(2.5,0)
\dashline{0.2}(4.5,0)(7.5,0)
\dashline{0.2}(9.5,0)(12,0)
\dashline{0.2}(3.5,1)(3.5,3)
\dashline{0.2}(8.5,1)(8.5,3)
\put(1.2,-1){$i$}\put(10.6,-1){$j$}\put(2.7,2){$m$}\put(8.7,2){$m'$}\put(5.8,-1){$k$}
\put(-.1,-.1){$\underbrace{\hspace*{2.7\unitlength}}$}
\put(9.4,-.1){$\underbrace{\hspace*{2.7\unitlength}}$}
\put(4.35,-.1){$\underbrace{\hspace*{3.2\unitlength}}$}
\end{picture}\\
\caption{The graph $P^{(m,m')}_{(i,k,j)}$.}\label{mm}
\end{center}
\end{figure}

\begin{lm}
  \label{rela} For $i,j,m,m'\geq 1$ and $k\geq 0$,
we have the following results on  the spectral radius of the tree $P^{(m,m')}_{(i,k,j)}$
(shown in Figure \ref{mm}).
\begin{enumerate}
\item
$
\lim\limits _{i,j
 \to \infty}\rho(P^{(m,m')}_{(i,k,j)})
\left\{
   \begin{array}{ll}
        >\frac{3}{2}\sqrt{2} & \mbox{ if } m,m'\ge 2 \mbox{ and } k\le m+m'; \mbox{ or one of } m \mbox{ and } m' \mbox{ is } 1,\\
    & \quad  k\leq m+m'-1 \mbox{ and } (m,m',k)\not=(1,1,1).\\
    =\frac{3}{2}\sqrt{2} & \mbox{ if } (m,m',k)=(1,1,1).\\
    <\frac{3}{2}\sqrt{2} & \mbox{ otherwise.}
   \end{array}
 \right.
$

\item
$
\lim\limits _{j \to \infty}\rho(P^{(m,m')}_{(m,k,j)})
\left\{
   \begin{array}{ll}
    >\frac{3}{2}\sqrt{2} & \mbox{ if } m\geq 2 \mbox{ and } k\leq m+m' -1,
(m,m',k)\not=(2,1,2), (2,2,3);\\
& \mbox{ or } m=1 \mbox {and } k\leq m'-2.\\
    <\frac{3}{2}\sqrt{2} & \mbox{ otherwise.}
   \end{array}
 \right.
$

\item
$ \rho(P^{(m,m')}_{(m,k,m')})
\left\{
   \begin{array}{ll}
    >\frac{3}{2}\sqrt{2} & \mbox{ if } m,m'\geq 2 \mbox{ and } k\leq m+m' -2,
(m,m',k)\not=(2,2,2);\\
& \mbox{ or } m=1 \mbox { and } k\leq m'-3.\\
    <\frac{3}{2}\sqrt{2} & \mbox{ otherwise.}
   \end{array}
 \right.
$
\end{enumerate}
\end{lm}

{\bf Proof } Similar to the computation in Lemma \ref{sptree}, we have
\begin{eqnarray*}
\phi_{P^{(m,m')}_{(i,k,j)}}
&=& (1, 1)A^iB_mA^kB_{m'}A^{j-1}
\left(\begin{array}[c]{c} p_1\\  q_1\end{array}\right)\\
&=& \frac{x_2^{i+j+1}}{x_2-x_1}(x_1^{2i}, 1)B_mA^kB_{m'}
\left(\begin{array}[c]{c} -x_1^{2(j+1)}\\ 1 \end{array}\right).
\end{eqnarray*}

By Lemma \ref{sqrt5},  the spectral radii  of all graphs considered in the lemma are  in $[2,\sqrt 5)$.
We can restrict $\lambda$ to this interval.

For item 1, let $\rho=\lim\limits _{i,j
 \to \infty}\rho(P^{(m,m')}_{(i,k,j)})$. Observe that  $\rho$ is the largest root of
the function
\begin{eqnarray*}
(0, 1)B_mA^kB_{m'}\left(\begin{array}[c]{c} 0\\ 1 \end{array}\right)
&=& \frac{1}{(x_2-x_1)^2} (-x_2\phi_{P_{m-1}}, d_{m}^{(2)})A^k
\left(\begin{array}[c]{c} x_1\phi_{P_{m'-1}}\\  d_{m'}^{(2)}\end{array}\right)\\
&=& \frac{x_2^{m+m'+k+2}}{(x_2-x_1)^4}f_{m,m',k}(\lambda).
\end{eqnarray*}
Here $f_{m,m',k}(\lambda)= (x_2^2-2 + x_1^{2m+2})
(x_2^2-2 + x_1^{2m'+2}) - x_1^{2k+2}(1-x_1^{2m})(1-x_1^{2m'})$.
The dominating term in $f_{m,m',k}(\lambda)$ is $(x_2^2-2)^2$.
We have $\lim\limits_{\lambda\to\infty}f_{m,m',k}(\lambda)=\infty$.

On one hand, to prove $\rho>\frac{3}{2}\sqrt{2}$,
we will show  $f_{m,m',k}(\frac{3}{2}\sqrt{2})< 0$.
On the other hand,  to prove $\rho<\frac{3}{2}\sqrt{2}$, we will show $f_{m,m',k}(\lambda)>0$ for all $\lambda\ge \frac{3}{2}\sqrt{2}$.

We assume $m,m'\geq 2$. Note that $x_2$ takes the value $\sqrt{2}$ at $\lambda=\frac{3}{2}\sqrt{2}$. If $k\leq m+m'$, then
\begin{eqnarray*}
  f_{m,m',k}\left(\frac{3}{2}\sqrt{2}\right)  &=& \frac{1}{2^{m+1}}\frac{1}{2^{m'+1}}
-\frac{1}{2^{k+1}}
\left (1-\frac{1}{2^m}\right)\left(1-\frac{1}{2^{m'}}\right)\\
&\leq & \frac{1}{2^{m+m'+1} }\left(\frac{1}{2} - \left (1-\frac{1}{2^m}\right)\left(1-\frac{1}{2^{m'}}\right) \right)\\
&<& 0,
\end{eqnarray*}
because of $(1-\frac{1}{2^{m}})(1-\frac{1}{2^{m'}})\geq \frac{3}{4}\cdot\frac{3}{4}>\frac{1}{2}$.
If $k\geq m+m'+1$, then  for $\lambda\geq \frac{3}{2}\sqrt{2}$
we have
$$f_{m,m',k}\left(\lambda\right)>x_1^{2(m+m')+4}-x_1^{2k+2}\geq 0.$$
Here we  applied  the fact $x_2^2\ge 2$ and $0<1-\frac{1}{2^m},1-\frac{1}{2^{m'}}<1$ for $\lambda\ge \frac{3}{2}\sqrt{2}$.

Now we assume one of $m$ and $m'$ is $1$, say $m=1$. If $k\le m+m'-1=m'$, then
$$
  f_{1,m',k}\left(\frac{3}{2}\sqrt{2}\right)=\frac{1}{2^{m'+3}}
-\frac{1}{2^{k+2}}\left(1-\frac{1}{2^{m'}}\right)\leq  \frac{1}{2^{m'+2} }\left(\frac{1}{2^{m'}}-\frac{1}{2} \right)\leq 0.
$$
The equality $f_{m,m',k}\left(\frac{3}{2}\sqrt{2}\right)=0$ holds if and only if $m=m'=k=1$.

If $k\ge m+m'=m'+1$, then for $\lambda\geq \frac{3}{2}\sqrt{2}$ we have
\begin{eqnarray*}
  f_{1,m',k}(\lambda) &=& (x_2^2-2 + x_1^{4})
(x_2^2-2 + x_1^{2m'+2}) - x_1^{2k+2}(1-x_1^{2})(1-x_1^{2m'})\\
&>& (x_2^2-2 + x_1^{4})x_1^{2m'+2} - x_1^{2m'+4}(1-x_1^{2})\\
&=&x_1^{2m'+2}(x_2^2-2)(1-x_1^4)\\
&\geq& 0.
\end{eqnarray*}

Overall, the proof of item 1 is completed.

For item 2, let $\rho'=\lim\limits _{j
 \to \infty}\rho(P^{(m,m')}_{(m,k,j)})$. A similar calculation shows that
 $\rho'$ is the largest root of
the following function
$$(x_1^{2m}, 1)B_mA^kB_{m'}  \left(\begin{array}[c]{c} 0\\ 1 \end{array}\right)
=\frac{(x_2^{m+1}-x_1^{m+1})x_2^{m'+k+1}}{(x_2-x_1)^4}g_{m,m',k}(\lambda),$$
 where $g_{m,m',k}(\lambda)=(x_2^2-2 + x_1^{2m})
(x_2^2-2 + x_1^{2m'+2}) - x_1^{2k+2}(1-x_1^{2m}(2-x_1^2))(1-x_1^{2m'})$.

For item 3, $\rho(H^{m,m'}_{m,k,m'})$ is the largest root of
the following function
$$(x_1^{2m}, 1)B_mA^kB_{m'}  \left(\begin{array}[c]{c} x_1^{2m'}\\ 1 \end{array}\right)
= \frac{(x_2^{m+1}-x_1^{m+1})(x_2^{m'+1}-x_1^{m'+1}) x_2^{k+1}}{(x_2-x_1)^4}h_{m,m',k}(\lambda),$$
where $h_{m,m',k}(\lambda)=(x_2^2-2 + x_1^{2m})(x_2^2-2 + x_1^{2m'}) -x_1^{2k+2}(1-x_1^{2m}(2-x_1^2))(1-x_1^{2m'}(2-x_1^2))$.

Items 2 and 3 can be proved by applying  similar arguments to  $g_{m,m',k}(\lambda)$ and
$h_{m,m',k}(\lambda)$ respectively.
The details are omitted here.
$\hfill\Box$

From Lemma \ref{rela} and \ref{sptree}, we get the following corollary. 
\begin{cor}\label{corela}
The following statements hold for $\rho_{m,k}$.
\begin{enumerate}
\item For $m\ge 2$, $\rho_{m,k}<\frac{3}{2}\sqrt 2$ if and only if $k\ge 2m+3$.
\item $\rho_{1,k}<\frac{3}{2}\sqrt 2$  if and only if $k\ge 4$.
\end{enumerate}
\end{cor}

\begin{cor} \label{cor42}
Suppose an open quipu $P_{(m_0,k_1,...,k_r,m_{r})}^{(m_0,...,m_r)}$ has spectral radius less than $\frac{3}{2}\sqrt{2}$. Then the following statements hold.
\begin{enumerate}
\item For $2\le i\le r-1$, we have $k_i\geq m_{i-1}+m_i$. Moreover if $m_{i-1},m_i\geq 2$, then  $k_i\geq m_{i-1}+m_i+1$.
\item  We have  $k_1 \geq m_{0}+m_1$ if $m_0\geq 2$;
and $k_1 \geq m_1-1$  if $m_0=1$.
\item We have $k_r \geq m_{r}+m_{r-1}$ if $m_r\geq 2$; and    $k_r \geq m_{r-1}-1$  if $m_r=1$.
\end{enumerate}
\end{cor}

The necessary conditions for $\rho(P_{(m_0,k_1,...,k_r,m_{r})}^{(m_0,...,m_r)})<\frac{3}{2}\sqrt{2}$ are quite good as evidenced by
the following theorem.

\begin{thm}\label{opqu2} Suppose that an open quipu $P_{(m_0,k_1,...,k_r,m_{r})}^{(m_0,...,m_r)}$  satisfies
\begin{enumerate}
\item $m_0,m_r\ge 2$;
\item $k_i\ge m_{i-1}+m_i+3$ for $2\le i\le r-1$;
\item $k_j\ge m_{j-1}+m_j+1$ for $j=1,r$.
\end{enumerate}
Then we have $\rho(P_{(m_0,k_1,...,k_r,m_r)}^{(m_0,...,m_r)})< \frac{3}{2}\sqrt{2}$.
\end{thm}

\noindent{\bf Proof } Denote $l_1=m_0+m_1+1$, $l_r=m_{r-1}+m_r+1$ and $l_i=m_{i-1}+m_i+3$ for $2\le i\le r-1$.
Let $G=P_{(m_0,l_1,...,l_r,m_r)}^{(m_0,...,m_r)}$.
By Lemma \ref{lm2.4}, we get $\rho(P_{(m_0,k_1,...,k_r,m_r)}^{(m_0,...,m_r)})\le \rho(G)$.
We have
\begin{eqnarray*}
\phi_{G}
&=& (1, 1)A^{m_0}B_{m_0}A^{l_1}B_{m_1}A^{l_2}\cdots A^{l_{r-1}}B_{m_{r-1}}A^{l_r}B_{m_r}A^{m_r-1}
\left(\begin{array}[c]{c} p_1\\ q_1 \end{array}\right)\\
&=& \frac{1}{(x_2-x_1)^3} (d_{m_0}^{(1)}x_1^{m_0}-x_2^{m_0+1}\phi_{P_{m_0-1}},  x_1^{m_0+1}\phi_{P_{m_0-1}}+d_{m_0}^{(2)}x_2^{m_0})\\
& &A^{l_1}B_{m_1}A^{l_2}\cdots A^{l_{r-1}}B_{m_{r-1}}A^{l_r+1}
\left( \begin{array}{ll} x_2^{m_r+1}\phi_{P_{m_r-1}}-d_{m_r}^{(1)}x_1^{m_r}\\ x_1^{m_r+1}\phi_{P_{m_r-1}}+d_{m_r}^{(2)}x_2^{m_r} \end{array}\right)\\
&=& \frac{(x_2^{m_0+1}-x_1^{m_0+1})(x_2^{m_r+1}-x_1^{m_r+1})}{(x_2-x_1)^3} (-d_{m_0-1}^{(1)},  d_{m_0-1}^{(2)})A^{l_1}B_{m_1}A^{l_2}\cdots A^{l_{r-1}}B_{m_{r-1}}A^{l_r+1}
\left( \begin{array}{ll} d_{m_r-1}^{(1)}\\ d_{m_r-1}^{(2)} \end{array}\right)\\
&=& \frac{(x_2^{m_0+1}-x_1^{m_0+1})(x_2^{m_r+1}-x_1^{m_r+1})}{(x_2-x_1)^3} (-d_{m_0-1}^{(1)}x_1^{m_0},  d_{m_0-1}^{(2)}x_2^{m_0})\\
&&(A^{m_1+1}B_{m_1}A^{m_1+2})\cdots
(A^{m_{r-1}+1}B_{m_{r-1}}A^{m_{r-1}+2})
\left( \begin{array}{ll} d_{m_r-1}^{(1)}x_1^{m_r}\\ d_{m_r-1}^{(2)} x_2^{m_r}\end{array}\right).
\end{eqnarray*}

Since $m_0,m_r\geq 2$, by Corollary \ref{corela},
we get
$d_{m_0-1}^{(2)}x_2^{m_0}> d_{m_0-1}^{(1)}x_1^{m_0}$  and
$d_{m_r-1}^{(2)}x_2^{m_r}> d_{m_r-1}^{(1)}x_1^{m_r}$
for all $\lambda\ge \frac{3}{2}\sqrt 2$.

Observe that  $A^{m+1}B_mA^{m+2}$ ($m\ge 1$) maps the region $\{(x,y)^T \mid x<y\}$ to itself. Now repeatedly apply this fact for $m=m_{r-1}, \ldots, m_1$.
We  get $\phi_{G}(\lambda)\ge 0$ for all $\lambda> \frac{3}{2}\sqrt 2$.
Thus, $$\rho(P_{(m_0,k_1,...,k_r,m_r)}^{(m_0,...,m_r)})\le\rho(G)< \frac{3}{2}\sqrt 2.$$

The proof is completed.
$\hfill\Box$

\subsection{Proofs of Theorems \ref{thm1} and  \ref{thm2}}
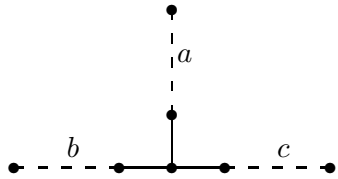
\begin{figure}[h]
\begin{center}
\setlength{\unitlength}{.7cm}
\begin{picture}(6,3)
\multiput(0,0)(2,0){2}{\circle*{0.18}}
\multiput(4,0)(2,0){2}{\circle*{0.18}}
\multiput(3,1)(0,2){2}{\circle*{0.18}}
\put(3,0){\circle*{0.18}}
\put(2,0){\line(1,0){2}}
\put(3,0){\line(0,1){1}}
\dashline{0.2}(0,0)(2,0)
\dashline{0.2}(4,0)(6,0)
\dashline{0.2}(3,1)(3,3)
\put(1,.2){$b$}\put(5,.2){$c$}\put(3.1,2){$a$}
\end{picture}
\end{center}
\caption{The T-shape trees}\label{tt}
\end{figure}

{\bf Proof of Theorem \ref{thm1}.}
Note that all the T-shape trees (see Figure \ref{tt}) have spectral radii less than $\frac{3\sqrt 2}{2}$ and satisfy  $3D> 2n-4$.

Now assume $r\ge 1$. Write $T$ as $P_{(k_0,k_1,\ldots, k_r,k_{r+1})}^{(m_0,m_1,\ldots, m_r)}$
with $m_i\geq 1$ for $i=0,1,\cdots,r$, $k_0\geq m_0$, $k_{r+1}\geq m_r$, $k_j\geq 0$ for $j=1,\ldots, r$.

\noindent
{\bf Case 1:} $r=1$. Here $T=P_{(k_0,k_1,k_2)}^{(m_0, m_1)}$. Without loss of generality,
we assume $m_0\leq m_1$.

Note that $T$ contains the
subgraph $P_{(m_0,k_1,m_1)}^{(m_0, m_1)}$. Since
$\rho(P_{(m_0,k_1,m_1)}^{(m_0,m_1)})\leq\rho(T)<\frac{3}{2}\sqrt{2}$,
by Item 3 of Lemma \ref{rela}, we must have $k_1\geq m_0+m_1-3$;  the equality 
holds  if and only if $m_0=1$.
Together with $k_0\geq m_0$ and $k_2\geq m_1$,  we get
\begin{eqnarray*}
  3D-(2n-4) &\geq& 3(k_0+k_1+k_2+1) -(2(k_0+k_1+k_2+m_0+m_1+2)-4)\\
&=& k_0+k_1+k_2-2m_0-2m_1+3\\
&\geq& m_0+ (m_0+m_1-3)+m_1 -2m_0-2m_1+3\\
&=&0.
\end{eqnarray*}
When the equality $  3D-(2n-4) =0$ holds,  we must have
$k_0=m_0=1$, $k_1=m_1-2$, and $k_2=m_1$. In this case, we get the
graph $P_{(1,m_1-2,m_1)}^{(1, m_1)}$, whose spectral radius is less than $\frac{3}{2}\sqrt 2$.

\noindent
{\bf Case 2:} $r=2$. Here $T=P_{(k_0,k_1,k_2,k_3)}^{(m_0, m_1, m_2)}$. Assume $m_0\leq m_2$ without loss of generality.

Note that $T$ contains the subgraph $P_{(m_0,k_1,k_2,m_2)}^{(m_0, m_1,
  m_2)}$. For any $i$ we have
$$\frac{3}{2}\sqrt{2}>\rho(T)\ge\rho(P_{(m_0,k_1,k_2,m_2)}^{(m_0,
  m_1, m_2)})>\lim_{i\to\infty}\rho(P_{(m_0,k_1,i)}^{(m_0, m_1)}).
$$
By Item 2 of Corollary \ref{cor42}, we have $k_1 \geq m_0+m_1-2$
with the equality if and only  if  $m_0=1$.
By symmetry, we also have $k_2 \geq m_1+m_2-2$
 with the equality  if and only if $m_2=1$.
Thus, we get
\begin{eqnarray*}
3D-(2n-4)&\geq& 3(k_0+k_1+k_2+k_3+2)-2(m_0+m_1+m_2+k_0+k_1+k_2+k_3+3)+4\\
&=&k_0+k_1+k_2+k_3-2(m_0+m_1+m_2)+4\\
&\geq& m_0+ (m_0+m_1-2)+(m_1+m_2-2)-2(m_0+m_1+m_2)+4\\
&=& 0.
\end{eqnarray*}
When the equality holds, we must have $k_3=m_2=k_0=m_0=1$ and
$k_1=k_2=m_1-1$.  We get the graph $T=P_{(1,m_1-1,m_1-1,1)}^{(1, m_1,
  1)}$, which has spectral radius greater than $\frac{3}{2}\sqrt 2$
shown as follows. For $m_1=1$, we can get $\rho(P_{(1,0,0,1)}^{(1, 1,
  1)})>\frac{3\sqrt 2}{2}$ by a straight calculation.  For $m_1\ge 2$,
by Lemma \ref{lm2.3}, and Item 3 of Lemma \ref{rela}, we have
$$\rho(P_{(1,m_1-1,m_1-1,1)}^{(1, m_1, 1)})=\rho(P_{(1,m_1-2,m_1+1)}^{(1,m_1+1)})>\frac{3}{2}\sqrt 2.$$

\noindent
{\bf Case 3:} $r\ge 3$.  Here $T=P_{(k_0,k_1,\ldots, k_r, k_{r+1})}^{(m_0, m_1,\ldots, m_r)}$.
Since $\rho(T)<\frac{3}{2}\sqrt{2}$, by Item 1 of Corollary \ref{cor42}, we must have
\begin{equation}
  \label{eq:20}
  m_{l-1}+m_l\leq k_l.
\end{equation}
By Items 2 and 3 of Corollary \ref{cor42}, we have
\begin{eqnarray}
  \label{eq:21}
  m_{0}+m_1&\leq& k_1+2, \\
  \label{eq:22}
  m_{r-1}+m_r&\leq& k_r+2.
\end{eqnarray}
Recall $m_0\leq k_0$ and $m_r\leq k_{r+1}$.
Summing up these two inequalities and
equations (\ref{eq:20}) (for $2\leq l \leq r-1$),
(\ref{eq:21}), (\ref{eq:22}),
 we get
\begin{equation}
  \label{eq:23}
  2 \sum_{l=0}^r m_l \leq  \sum_{i=0}^{r+1} k_i+4.
\end{equation}
Hence, we have
\begin{eqnarray*}
3D-(2n-4)&\geq& 3\left(\sum_{i=0}^{r+1}k_i+r\right)-2\left(\sum_{j=0}^rm_j + \sum_{i=0}^{r+1}k_i +r +1 \right)+4\\
&=&\sum_{i=0}^{r+1}k_i -2 \sum_{j=0}^rm_j  + r+2\\
&\geq& r-2 \\
&>&0.
\end{eqnarray*}
The proof of Theorem \ref{thm1} is completed.
\hfill$\Box$

{\bf Proof of Theorem \ref{thm2}.}
Let $L=C_{(k_1,...,k_r)}^{(m_1,...,m_r)}$, where $k_i\ge 0$ and $m_i\ge 1$ for $i=1,...,r$.
For convenience, we write $m_0=m_r$.

First, we prove the lower bound of $D(L)$. Denote $m=\max \{ m_1,...,m_r\}$.
We have
\begin{eqnarray}
\label{order}
n&=&r+\sum_{i=1}^r m_i+\sum_{i=1}^r k_i,\\
\label{diam}
D&\ge& m+\left\lfloor \frac{1}{2}(r+\sum_{i=1}^r k_i)\right\rfloor>\frac{1}{2}(r+\sum_{i=1}^r k_i).
\end{eqnarray}

 By the condition $\rho(L)\le \frac{3}{2}\sqrt{2}$ and  Item 1 of Corollary \ref{cor42},
we have $m_{i-1}+m_i\leq k_i$
for all $1\leq i\leq r$. We get
 \begin{equation}\label{remk}
 2\cdot\sum_{i=1}^r m_i\le \sum_{i=1}^r k_i.
 \end{equation}

Let $\bar m=\frac{\sum_{i=1}^r m_i}{r}$. Combining the inequalities (\ref{order}), (\ref{diam}), and (\ref{remk}), we get
$$\frac{n}{n-2D}>\frac{r+3\cdot\sum_{i=1}^r m_i}{\sum_{i=1}^r m_i}=3+\frac{1}{\bar m}.$$
Solving for $D$, we get 
$$D>\frac{2+\frac{1}{\bar m}}{2\cdot(3+\frac{1}{\bar m})}n>\frac{n}{3}.$$

Now we prove the upper bound $\frac{2n-4}{3}$ for  $D(L)$.

If $r=1$, then $L=C_{(k)}^{(m)}$. We have $\rho(L)=\rho_{m,k}$.
By Corollary \ref{corela} and $n\ge 13$, we have $k\geq 2m+3$ and
\begin{eqnarray*}
  3D-2n+4 &=& 3\left(m+\left\lfloor \frac{k+1}{2}\right\rfloor \right)
-2(m+k+1)+4\\
&=&m+3\left\lfloor \frac{k+1}{2}\right\rfloor -2k+2.
\end{eqnarray*}
When $k=2t$ even, since $2t=k\geq 2m+3$, we get
$$3D-2n+4 =m-t+2\leq 0.$$
When $k=2t+1$ odd, since $2t+1\geq 2m+3$,
we get
$$3D-2n+4 =m-t+3\leq 2.$$

 Here we get two exception cases to $3D\leq 2n-4$:
$k=2m+3$ and $k=2m+5$ (the graphs are shown in Figure \ref{excp}).

Now we consider the case $r\geq 2$. Let $m$ (or $m'$) be the first (or the second) 
largest number in $\{m_1,\ldots, m_r\}$ respectively. Let $L'$ be the graph obtained
from $L$ by removing all pendent paths other than the two longest ones.
Let $g$ denote the length of the unique cycle in $L$.  Let
$L''=C_{(\lfloor\frac{g}{2}\rfloor-1,\lceil\frac{g}{2}\rceil-1)}^{(m,m')}$.
By Lemma \ref{lm2.1} and Lemma \ref{csh}, we have
$$\rho(L)\geq \rho(L')\geq \rho(L'').$$
We observe that $n(L)\geq n(L'')$ and $D(L)\leq D(L'')$. It suffices to
prove $3D(L'')\leq 2n(L'')-4.$

If $g=2k$ even, then
$\rho(L'')=\rho(C^{(m,m')}_{(k-1,k-1)} )=
\rho(P_{(m+1,k-3,m'+1}^{m+1,m'+1}).$
Since $\rho(P_{(m+1,k-3,m'+1}^{m+1,m'+1})=\rho(L'')< \frac{3}{2}\sqrt{2}$,
by Item 3 of Lemma \ref{rela}, we get
$m+m'+2\le  k-2$ unless $m=m'=1$ and $k=5$.
We will consider the special case later. For general case,  we have
\[3D(L'')-2n(L'')+4 = 3(m+m'+k)-2(m+m'+2k) +4 = m+m'-k +4
 \leq 0. \]
When  $m=m'=1$ and $k=5$, we have $L''=C_{(4,4)}^{(1,1)}$
and $3D(L'')= 2n(L'')-3.$  Since  $n(L)\geq 13>12=n(L'')$, we have
\[3D(L)-2n(L)+4< 3D(L'')-2n(L'')+4 =1.\]

When $g=2k+1$ is odd, let $L'''=C^{(m-1, m')}_{(k,k)}$.
Since $L'''$ can be obtained from $L''$ by deleting a leaf vertex and
subdividing an internal edge, by Lemma \ref{lm2.4}, we have $\rho(L''')<\rho(L'')$.
We also observe that $n(L''')=n(L'')$ and $D(L''')=D(L'')$.
By the previous cases, we have $3D(L''')\leq 2n(L''')-4$.
Thus  $3D(L'')\leq 2n(L'')-4.$ We are done. \hfill $\square$

\section{Application to diameter $\frac{n}{2}\le D\le\frac{2n-4}{3}$}
We have the following lemma.
\begin{lm}\label{roro}
For $m\geq 1$, let $\rho_m=\lim_{k\to\infty} \rho_{m,k}$.  We have
 $\rho_{m+1}>\rho_{m,k}$ holds for  $k\ge 2m+5$.
\end{lm}

\noindent {\bf Proof } Recall $\rho_{m,k}$ is the largest root of
$\frac{d_m^{(2)}}{d_m^{(1)}}=x_1^{k-1}$.  Thus $\rho_{m+1}$ is the
largest roots of $d_{m+1}^{(2)}=0$ while $\rho_{m,2m+5}$ is the
largest roots of $\frac{d_m^{(2)}}{d_m^{(1)}}=x_1^{2m+4}$. 
 Let $f=f(\lambda)$ be a function of $\lambda$. The
notation $f|_{\lambda_0}$ means the value of $f$ at
$\lambda=\lambda_0$.  We have
\begin{equation}
\label{eq31}
\left .(x_2^2-2)+x_1^{2m+4} \right |_{\rho_{m+1}}=0
\end{equation}
and
\begin{equation}
\label{eq32}
\left .(x_2^2-2)+x_1^{2m+4}\left( x_2^2-1+x_1^{2m+2}(2-x_1^2)\right)\right |_{\rho_{m,2m+5}}=0.
\end{equation}

We get
\begin{eqnarray*}
\left .(x_2^2-2)+x_1^{2m+4} \right |_{\rho_{m,2m+5}}&=& \left .x_1^{2m+4}- x_1^{2m+4}\left( x_2^2-1+x_1^{2m+2}(2-x_1^2)\right )\right|_{\rho_{m,2m+5}}\\
&=& \left .x_1^{2m+4}\left( 2- x_2^2-x_1^{2m+2}(2-x_1^2)\right )\right |_{\rho_{m,2m+5}}\\
&=& \left .x_1^{2m+4}\left[ x_1^{2m+4}\left( x_2^2-1+x_1^{2m+2}(2-x_1^2)\right)-x_1^{2m+2}(2-x_1^2)\right ]\right |_{\rho_{m,2m+5}}\\
&=& \left .x_1^{4m+6}\left(x_1^{2m+4}(2-x_1^2)-1\right )\right |_{\rho_{m,2m+5}}\\
&\le & \left .x_1^{4m+6}\left(x_1^6(2-x_1^2)-1\right )\right |_{\rho_{m,2m+5}}\\
&<& \left .x_1^{4m+6}\left(\frac{3x_1^6}{2}-1\right )\right |_{\rho_{m,2m+5}}\\
&<&0.
\end{eqnarray*}

In the last step we use $x_1^2|_{\rho_{m,2m+5}}<\frac{\sqrt 5-1}{2}$,
since $\rho_{m,2m+5}>\sqrt{2+\sqrt 5}$.
Thus, we have $\rho_{m+1}>\rho_{m,2m+5}\ge\rho_{m,k}$. The proof is completed.
$\hfill\Box$

Lemma 4.3 of \cite{CDK} can be generalized to the following lemma.
The proof is similar and will be omitted.
\begin{lm} \label{lmlast}
If a minimizer graph with $n$ vertices and
diameter $D$ with $n\ge D+2$ and $\frac{n}{2}\le D\le\frac{2n-4}{3}$ is a subgraph of an ($D-\lfloor\frac{n}{2}\rfloor$)-Urchin graph but not of an ($D-\lfloor\frac{n}{2}\rfloor$)-Laundry graph, then it is $C_{(n-D-1,n-D-1)}^{(D-\lfloor\frac{n}{2}\rfloor,D-\lceil\frac{n}{2}\rceil)}$.
\end{lm}

\noindent{\bf Proof  of Theorem \ref{thm3}.}
To apply Lemma \ref{lmlast}, it  suffices to prove the following two claims for $n\geq 13$.

{\bf Claim 1.} $G_{n,D}^{min}$ must be a closed quipu.

{\bf Claim 2.}  The longest pendent path of $G_{n,D}^{min}$  has length
at most $D-\lfloor\frac{n}{2}\rfloor$.

First we prove Claim 1. Consider the graph
$C_{(n-D-1,n-D-1)}^{(D-\lfloor\frac{n}{2}\rfloor,D-\lceil\frac{n}{2}\rceil)}$. Let
$m=D-\lfloor\frac{n}{2}\rfloor$ and
$k=n-D-1$. Since $n\ge 13$ and $\frac{n}{2}\le D<\frac{2n-4}{3}$, we
have $k\ge 4$ and $k>2m+3$.
By Corollary \ref{corela},  we have
$$\rho(G_{n,D}^{min})\le\rho\left(C_{(n-D-1,n-D-1)}^{(D-\lfloor\frac{n}{2}\rfloor,D-\lceil\frac{n}{2}\rceil)}\right)\le\rho\left(C_{(k,k)}^{(m,m)}\right)<\frac{3\sqrt 2}{2}.$$

So, $G_{n,D}^{min}$ is  either a dagger, an open quipu, or a closed quipu.
The minimizer graph $G_{n,D}^{min}$ can not be a dagger since  all daggers with $n\ge 6$
do not satisfy $3D\leq 2n-4$.
By Theorem \ref{thm1}, $G_{n,D}^{min}$ can not be an open quipu either. Hence, $G_{n,D}^{min}$  must be a closed quipu.

Now we prove Claim 2. Since $3D\le 2n-7$,  we have
 $n-D-1\ge 2(D-\lfloor\frac{n}{2}\rfloor)+5$.
Suppose that $G_{n,D}^{min}$ has a pendent path of length $m'>D-\lfloor\frac{n}{2}\rfloor$.
By Lemma \ref{roro}, we have $$\rho\left(C_{(n-D-1,n-D-1)}^{(D-\lfloor\frac{n}{2}\rfloor,D-\lceil\frac{n}{2}\rceil)}\right)<\rho_{m'}<\rho(G_{n,D}^{min}).$$ Contradiction!
The proof of two claims are finished. Applying Lemma \ref{lmlast}, we are done.
\hfill $\square$

\end{document}